\documentclass[12pt,a4paper]{article}

\makeatletter
\long\def\@makecaption#1#2{%
  \vskip\abovecaptionskip
  \sbox\@tempboxa{#1. #2}%
  \ifdim \wd\@tempboxa >\hsize
    #1. #2\par
  \else
    \global \@minipagefalse
    \hb@xt@\hsize{\hfil\box\@tempboxa\hfil}%
  \fi
  \vskip\belowcaptionskip}
\makeatother

\usepackage[english]{babel}
\usepackage{amssymb}
\usepackage{amsmath}
\usepackage[pdftex]{graphicx}
\usepackage{graphics}
\usepackage{times}
\usepackage{url}
\usepackage{multicol}

\newtheorem{tm}{Theorem}
\newtheorem{defi}{Definition}
\newtheorem{lem}{Lemma}
\newtheorem{cor}{Corollary}

\begin{document}

 \begin{center} {\large {\bf Circular Surfaces $\mathcal {CS}(\alpha,p)$}}

\bigskip

 SONJA GORJANC\\
{\small University of Zagreb, Faculty of Civil Engineering,\\ Ka\v{c}i\'{c}eva 26, 10000 Zagreb, Croatia\\
e-mail: sgorjanc@master.grad.hr}

\medskip

 EMA JURKIN\\
{\small University of Zagreb, Faculty of Mining, Geology and Petroleum Engineering,\\ Pierottijeva 6, 10000 Zagreb, Croatia\\
e-mail: ema.jurkin@rgn.hr}
\end{center}

\bigskip

\noindent {\small {\bf Abstract}.\\
In this paper we define and construct a new class of algebraic surfaces in three-dimensional Euclidean space generated by a curve   and a congruence of circles. We study their properties and visualize them with the program {\it Mathematica}.

\medskip
\medskip

}\\

\section{Introduction}

In  article \cite{Izumiya}  the definition of a {\it circular surface} is given as follows: 
\begin{defi}
A {\it circular surface} is (the image of) a map $V:I\times \mathbb R/2\pi \mathbb Z \longrightarrow \mathbb R^3$ defined by
\begin{equation}\label{cs}
V(t,\theta)=\gamma(t)+r(t)(\cos \theta\mathbf a_1(t)+\sin \theta\mathbf a_2(t)),
\end{equation}
where $\gamma, \mathbf a_1, \mathbf a_2 \colon I\longrightarrow \mathbb R^3$ and $r \colon I\longrightarrow \mathbb R_{>0}$.
\end{defi}
It is assumed that $\langle \mathbf a_1,\mathbf a_1\rangle =\langle \mathbf a_2,\mathbf a_2\rangle=1$ and $\langle \mathbf a_1,\mathbf a_2\rangle=0$ for all $t\in I$, where $\langle ,\rangle$ denotes the canonical inner product on $R^3$. The curve $\gamma$ is called   a {\it base curve} and the pair of curves $\mathbf a_1$, $\mathbf a_2$ is said to be   a {\it director frame}. The standard circles $\theta\mapsto r(t)(\cos\theta\mathbf a_1(t)+\cos\theta\mathbf a_2(t))$ are called {\it generating circles}. The base curve actually consists of the centers of the generating circles. Furthermore, in \cite{Izumiya} the authors classify and investigate the properties of this one-parameter family of standard circles with a fixed radius by using the methods of differential geometry.  

In this paper we construct one new class of circular surfaces.  First we define three types of congruences of circles and study their properties. Then, for a given congruence $\mathcal C(p)$  of a certain type and a given curve $\alpha$  we define a circular surface as the system of circles from $\mathcal C(p)$  that intersect $\alpha$. We obtain circular surfaces with  variable radii of the generating circles, and in that sense this paper is an extension of \cite{Izumiya}. For  algebraic curves $\alpha$ we investigate  algebraic properties of the circular surfaces -- the main  scientific contribution of this paper is stated in Theorem~\ref{tm1}. We use a relatively unusual method since we combine two different approaches: one is the analytical and the other is the synthetic approach. The analytical method is widely used and necessary when it comes to computer modeling and visualization (we use {\it Mathematica}). On the other hand, the synthetic method is often considered to be obsolete but in many cases it is useful and offers a short and  elegant solution to the problem.

\section{Congruence of circles $\mathcal C(p)$}

A congruence of curves is a set  $\mathcal C$ of curves in a three-dimensional space (projective, affine or Euclidean) depending on two parameters. According to \cite{Eisenhart}, it can be defined by equations of the form
\begin{equation}
\label{eisen}
x=f_1(t,u,v),\quad y=f_2(t,u,v),\quad z=f_3(t,u,v),
\end{equation}
where $f_1$, $f_2$, $f_3$ are the functions of class $C^1$ within a suitable region, $u$ and $v$ denote the parameters determining the curve, and $t$ is the parameter which determines the point on the curve.

\noindent In this paper we consider a {\it congruence of circles} $\mathcal C(p)$ that consists of circles in Euclidean space $\mathbb E^3$  passing through two given points $P_1$, $P_2$. The points $P_1$, $P_2$ lie on the axis $z$ and are given by the coordinates $(0,0,\pm p)$, where $p=\sqrt{q}$, $q\in\mathbb R$. If $q$ is greater, equal or less than zero, i.e. the points $P_1$, $P_2$ are real and different, coinciding (the axis $z$ is the tangent line of all circles of the congruence) or imaginary, the congruence $\mathcal C(p)$ is called an {\it elliptic},   {\it parabolic} or  {\it hyperbolic} congruence, respectively (see Fig. 1).

\begin{center}
\includegraphics[scale=0.75]{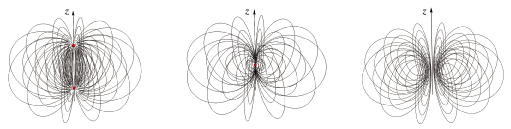}
\vspace{-0.5cm}
\begin{multicols}{3}
a

b

c
\end{multicols}
\end{center}
\vspace{-0.2cm}
Figure 1: 
\parbox[t]{120mm}{Some circles of an elliptic, parabolic or hyperbolic congruence $\mathcal C(p)$ are shown in a, b and c, respectively.}

\bigskip

\noindent Every plane $\zeta$ through the axis $z$ cuts the congruence $\mathcal C(p)$ into one pencil of circles through  the points $P_1$ and $P_2$ which is denoted $(P_1,P_2)$. This pencil is elliptic, parabolic or hyperbolic if $q$ is greater, equal or less than 0, respectively. In the plane $\zeta (\varphi)$, where $\varphi\in [0,2\pi)$ is the angle between the plane $\zeta$ and the axis $x$, we consider the Cartesian coordinate system $O(\rho,z)$ (see Fig. 2). 

\begin{center}
\includegraphics[scale=0.87]{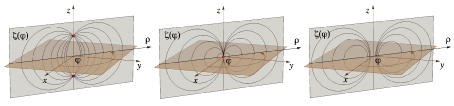}
\vspace{-0.5cm}
\begin{multicols}{3}
a

b

c
\end{multicols}
\vspace{-0.2cm}
Figure 2
\end{center}

\bigskip

\noindent For every circle $c\in (P_1,P_2)$ in the plane $\zeta(\varphi)$, the $\rho$-coordinate of its center $C$ is equal to $\sqrt{r^2-p^2}$, where $r$ is the radius of the circle $c$.

\noindent Thus, the congruence $\mathcal C(p)$ is defined by  the following equations:
\begin{align}
\label{congruence_xyz1}
&x(\theta,r,\varphi)=\cos\varphi(r \cos\theta+\sqrt{r^2-p^2}),\nonumber\\
&y(\theta,r,\varphi)=\sin\varphi(r \cos\theta+\sqrt{r^2-p^2}),\nonumber\\
&z(\theta,r,\varphi)=r \sin\theta,
\end{align}
where $\varphi,\theta\in [0,2\pi)$,  $r\in\mathbb R^+$ for $q\leqslant 0$, and $r\geqslant p$ for $q>0$.

\bigskip

\noindent It follows from eq. (\ref{eisen}) that the curves of a congruence are the intersections of  two systems of surfaces when $u$ and $v$ respectively are constants. We derive eqs. (\ref {congruence_xyz1})  on the basis that the circles of $\mathcal C(p)$ for constant $\varphi$ lie in the planes $\zeta(\varphi)$ (see Fig. 2). On the other hand, it is clear from eqs. (\ref {congruence_xyz1}) that for a constant $r$ the circles of $\mathcal C(p)$ lie on a torus $\tau (r)$. For  a variable $r$, the system of $\tau(r)$ contains  the spindle, horn or ring tori if $\mathcal C(p)$ is an elliptic, parabolic or hyperbolic congruence, respectively (see Fig. 3).

\begin{center}
\includegraphics[scale=0.8]{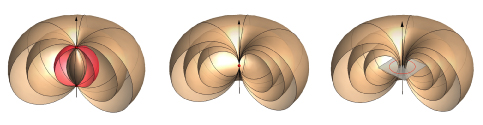}
\vspace{-0.5cm}
\begin{multicols}{3}
a

b

c
\end{multicols}
\end{center}
\vspace{-0.2cm}
Figure 3: 
\parbox[t]{120mm}{Some tori  of the system $\tau(r)$, for $q$ greater, equal or less than $0$, are shown in  a, b and c, respectively.
For $q>0$, $\tau(p)$ is a double sphere through the points $P_1$ and $P_2$ (a). For $q<0$, $\tau(0)$ is a ruled surface with the isotropic rulings that intersect on one real circle $c(0)$ given by the equations $x^2+y^2=-p^2$ and $z=0$ (c). }

\bigskip

\noindent For every point $A$ with cylindrical coordinates $(\rho_0,\varphi_0,z_0)$, $\rho_0\neq 0$, there exists a unique circle $c^A(p)\in\mathcal C(p)$ passing through the points $A$, $P_1$ and $P_2$. The $\rho$-coordinate of the center $C$ of this circle  $c^A(p)$ is $\rho_C^A=\frac{\rho_0^2+z_0^2-p^2}{2\rho_0}$.
 Therefore, if we convert the cylindrical coordinates to  Cartesian coordinates, then  $\forall A=(x_0,y_0,z_0)$, where $(x_0,y_0)\neq (0,0)$, the parametric equations of the circle $c^A\in \mathcal C(p)$ are the following:

{\small
\begin{align}
\label{Cartesian_cA}
&x(\theta)=\frac{x_0}{2(x_0^2+y_0^2)}\Big(\sqrt{4p^2(x_0^2+y_0^2)+(x_0^2+y_0^2+z_0^2-p^2)^2}\cos\theta+x_0^2+y_0^2+z_0^2-p^2\Big)\nonumber\\
&y(\theta)=\frac{y_0}{2(x_0^2+y_0^2)}\Big(\sqrt{4p^2(x_0^2+y_0^2)+(x_0^2+y_0^2+z_0^2-p^2)^2}\cos\theta+x_0^2+y_0^2+z_0^2-p^2\Big)\nonumber\\
&z(\theta)= \sqrt{p^2+\frac{(x_0^2+y_0^2+z_0^2-p^2)^2}{4(x_0^2+y_0^2)}}\sin\theta, \quad\quad \theta\in [0,2\pi).
\end{align}
}
\medskip

\noindent For every point at infinity $A^\infty$ $(A^\infty\notin z)$, the circle $c^{A^\infty}(p)$ splits into the axis $z$ and the line at infinity in the plane $\zeta$ through the axis $z$ and the point $A^\infty$.

\bigskip

\noindent {\bf Singular points of $\mathcal C(p)$}

\noindent A point is the {\it singular point} of a congruence if infinitely many curves pass through it.  The singular points of $\mathcal C(p)$ are the points on the axis $z$ and the absolute points of $\mathbb E^3$. Namely,
\begin{itemize}
\item
for the points $P_1,P_2\in z$, the circles $c_{P_i}$ form the whole congruence $\mathcal C(p)$;
\item if $Z\in z$, $Z\neq P_i$, $c^{Z}(p)$ form a singly infinite system of splitting  circles. Every plane $\zeta$ cuts this system into  $z$ and the line at infinity;
\item if $A_1,A_2$ are  the absolute points in the plane $\zeta$, $c^{A_1,A_2}(p)$ form a pencil of circles in the plane $\zeta$.
\end{itemize}

\bigskip

\noindent {\bf  $\mathcal C(p)$ is a normal curve congruence}

\noindent If all curves of a congruence are orthogonal to a singly infinite system of surfaces, the congruence is said to be {\it normal} and its equations satisfy the differential conditions given in \cite[p. 123]{Salmon2}. From eqs. (\ref{Cartesian_cA}) it is not hard to check  that these conditions are fulfilled for the congruence $\mathcal C(p)$.  
\noindent Even easier, according to the properties of {\it Apollonian circles} (\cite{wiki_Apollonian}), it is clear that every circle of $\mathcal C(p)$ cuts orthogonally one singly infinite system of spheres with the centers on the axis $z$ (see Fig. 4). 

\begin{center}
\includegraphics[scale=0.7]{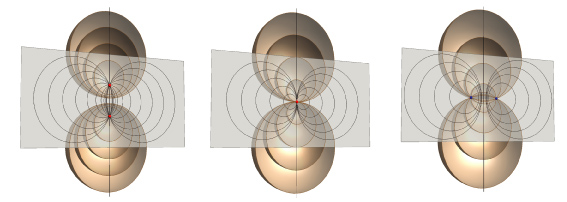}
\end{center}

\noindent Figure 4: \parbox[t]{120mm}{The figure illustrates the following property: every circle of an elliptic, parabolic or hyperbolic congruence $\mathcal C(p)$ is orthogonal to one hyperbolic, parabolic or elliptic pencil of spheres, respectively.}

\section{Circular surface $\mathcal {CS}(\alpha,p)$}

\begin{defi}
\label{def2}
Let $\alpha \colon I\rightarrow\mathbb R^3$, $I\subseteq \mathbb R$, be a piecewise-differentiable curve in $\mathbb E^3$ that does not contain any singular point of $\mathcal C(p)$. The circular surface $\mathcal C(\alpha,p)$ is a one-parametric system of  circles $c \in \mathcal C(p)$ that intersect the curve $\alpha$.
\end{defi}

\begin{center}
\includegraphics[scale=0.8]{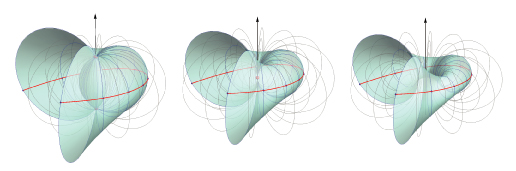}
\vspace{-0.5cm}
\begin{multicols}{3}
a

b

c
\end{multicols}
\end{center}
\vspace{-0.2cm}
Figure 5: 
\parbox[t]{120mm}{The figure illustrates Definition \ref{def2}. Through every point on the curve $\alpha$ one and only one circle of the congruence passes.}

\bigskip

\noindent If $\alpha$ contains any singular point of $\mathcal C(p)$, the one-parametric system of  circles $c \in \mathcal C(p)$ that intersect the curve $\alpha$ contains circles that have to be excluded from our further consideration. For example, if $\alpha$ passes through the point $P_i$, this system is the two-parametric set $\mathcal C(p)$. Thus, let us define a one-parametric set $\mathcal K(\alpha,p)$ that  further in the text we will consider as the complete system of circles $c\in \mathcal C(p)$ that cut $\alpha$.

\noindent If $\alpha$ is a piecewise-differentiable curve in $\mathbb E^3$, $\mathcal K(\alpha,p)$ consists of all circles $c\in\mathcal C(p)$ that cut $\alpha$ at the points that are regular points of $\mathcal C(p)$, and the following lines:
\begin{itemize}
\item[-] if $\alpha$ contains a point $Z\in z$, $Z\neq P_i$, $\mathcal K(\alpha,p)$ contains the pencil of lines through the point $(0\!:\!0\!:\!1\!:0)$ in the plane at infinity $\omega$. The line at infinity that lies in the plane $\zeta$ determined by the tangent line of $\alpha$ at $Z$ is denoted by $l_Z$;
\item[-] if $\alpha$ passes through a pair of  absolute points $(A_1,A_2)$, $\mathcal K(\alpha,p)$ contains all pairs of isotropic lines that lie in  isotropic planes $\iota_{1,2}$ through $z$ and $(A_1,A_2)$;

\item[-] if $\alpha$ passes through the point $P_i$, $\mathcal K(\alpha,p)$ contains the pencil of lines at the plane at infinity $\omega$ and the pencil of circles $c\in \mathcal C(p)$  in the plane $\pi$ that passes through the axes $z$ and the tangent line of $\alpha$ at $P_i$. The circle of this pencil that has the same tangent line at $P$ as the curve $\alpha$ is denoted by   $c_P$.
\end{itemize}

In the following notation, $\mathcal K(\alpha,p)$, $\omega$, $\iota_{1,2}$ and $\pi$ are implied as the sets of lines.

\begin{defi}
\label{def3}
Let $\alpha$ be a piecewise-differentiable curve in $\mathbb E^3$ that contains singular points of $\mathcal C(p)$. Then,
\begin{itemize}
\item[-] if $\alpha$ contains a point $Z$ on the axis $z$, $Z\neq P_i$, the circular surface  $\mathcal {CS}(\alpha,p)=(\mathcal K(\alpha,p) \,\,\backslash\,\, \omega)\cup\{l_Z\}$;
\item[-] if $\alpha$ contains the pair of absolute points $(A_1,A_2)$, the circular surface  $\mathcal {CS}(\alpha,p)=\mathcal K(\alpha,p) \,\,\backslash\,\, \iota_{1,2}$;
\item[-] if $\alpha$ contains the point $P_i$,  $\mathcal {CS}(\alpha,p)=(\mathcal K(\alpha,p) \,\,\backslash\,\, (\pi \cup \omega))\cup\{c_P\}$.
\end{itemize}
\end{defi}

\subsection{Parametric equations of $\mathcal {CS}(\alpha,p)$}

\noindent Let  the curve $\alpha\colon I\rightarrow\mathbb R^3$, $I\subseteq \mathbb R$, be given by
$\alpha(t)=(\alpha_1(t), \alpha_2(t), \alpha_3(t))$, where $\alpha_i \in C^1(I)$, and let its  normal projection  on the plane $z=0$ be denoted by $\alpha_{xy}(t)=(\alpha_1(t), \alpha_2(t), 0)$.

\noindent If we substitute $(\alpha_1(t), \alpha_2(t), \alpha_3(t))$ for $(x_0,y_0,z_0)$ in (\ref{Cartesian_cA}), then after some calculations  the following  parametric equations of  $\mathcal {CS}(\alpha,p)$ hold: 
\begin{align}
\label{par.surface1}
x(t,\theta)&= \frac{\alpha_1(t)}{2\left \| \alpha_{xy}(t) \right \|^2}\left(\sqrt{4p^2\left \| \alpha_{xy}(t) \right \|^2 +(\left \| \alpha(t) \right \| ^2-p^2)^2}\cos\theta+\left \| \alpha(t) \right \| ^2-p^2\right)\nonumber\\
y(t,\theta)&= \frac{\alpha_2(t)}{2\left \| \alpha_{xy}(t) \right \|^2}\left(\sqrt{4p^2\left \| \alpha_{xy}(t) \right \|^2 +(\left \| \alpha(t) \right \| ^2-p^2)^2}\cos\theta+\left \| \alpha(t) \right \| ^2-p^2\right)\nonumber\\
z(t,\theta)&= \sqrt{p^2+\frac{(\left \| \alpha(t) \right \| -p^2)^2}{4\left \| \alpha_{xy}(t) \right \| ^2}}\sin\theta, \quad\quad (t, \theta) \in I \times [0,2\pi).
\end{align}
The equations above can be written in the following form:
\begin{equation}\label{par.surface2} \mathcal{CS}^{\alpha\,p}(t,\theta)=\gamma^{\alpha\,p}(t)+r^{\alpha\,p}(t)(\cos \theta\mathbf a_1^{\alpha\,p}(t)+\sin \theta\mathbf a_2^{\alpha\,p}(t)), \end{equation} 
where
\begin{align}
& \gamma^{\alpha\,p}(t)=\frac{\lVert\alpha(t)\rVert^2-p^2}{2\lVert\alpha_{xy}(t)\rVert^2}\,\alpha_{xy}(t), 
\quad \quad \quad \mathbf a_1^{\alpha\,p}(t)=\frac{\alpha_{xy}(t)}{\lVert\alpha_{xy}(t)\rVert},\nonumber \\
& r^{\alpha\,p}(t)=\sqrt{p^2+\frac{(\lVert\alpha(t)\rVert^2-p^2)^2}{4\lVert\alpha_{xy}(t)\rVert^2}},
\quad \mathbf a_2^{\alpha\,p}(t)=(0,0,1)=\mathbf  k. \nonumber 
\end{align}
The equation (\ref{par.surface2}) is obviously of the form  (\ref{cs}) and therefore the  surface $\mathcal {CS}(\alpha,p)$  belongs to the class of {\it circular surfaces}.

\bigskip

\noindent The equations (\ref{par.surface1}) or (\ref{par.surface2}) allow numerous visualizations of the circular surfaces $\mathcal{CS}(\alpha,p)$ in the program {\it Mathematica}. Some examples are shown in  Figures 6 and 7.

\begin{center}
\includegraphics[scale=0.98]{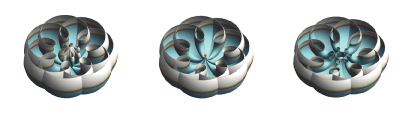}
\vspace{-0.5cm}
\begin{multicols}{3}
a

b

c
\end{multicols}
\end{center}
\vspace{-0.2cm}
Figure 6: \parbox[t]{120mm}{Parts of $\mathcal{CS}(\alpha,p)$ for $\alpha(t)=(9\cos t-4\cos \frac{9t}{2},9\sin t-4\cos \frac{9t}{2},0)$ and $p=4,0,4i$ are shown in  a, b and c, respectively.}

\begin{center}
\includegraphics[scale=0.98]{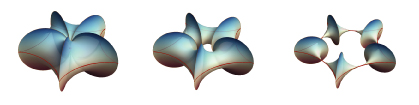}
\vspace{-0.5cm}
\begin{multicols}{3}
a

b

c
\end{multicols}
\end{center}
\begin{center}
\includegraphics[scale=0.98]{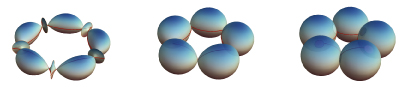}
\vspace{-0.5cm}
\begin{multicols}{3}
d

e

f
\end{multicols}
\end{center}
\vspace{-0.2cm}
Figure 7: \parbox[t]{120mm}{Six surfaces $\mathcal{CS}(\alpha,p)$ for $\alpha(t)=(4\cos t+\cos 4t,4\sin t-\sin 4t,0)$ and $p=0,2i,3i,4i,5i,6i$ are shown in  a, b, c, d, e and f, respectively.}

\subsection{Properties of algebraic surfaces $\mathcal {CS}(\alpha,p)$}

\bigskip
\begin{lem}\label{L1}
If  $\mathbf l$ is a straight line that does not intersect the axis $z$, $\mathcal{CS}(\mathbf l,p)$ is a parabolic cyclide through the axis $z$.
\end{lem}
{\sc Proof.} Without loss of generality we can assume that the line $\mathbf l$ is given by $$\mathbf l (t)=(1,t,bt+c), \quad t, b, c \in\mathbb R.$$ 
It follows from (\ref{par.surface1}) that the parametric equations of $\mathcal{CS}(\mathbf l,p)$ are
\begin{align}
\label{par.par.cyc}
x(t,\theta)&= \frac{1 - p^2 + t^2 + \left(c + b t\right)^2+\sqrt{4\left (1 + t^2\right)p^2 +\left(1 - p^2 + t^2 + \left(c + b t\right)^2\right)^2}\cos\theta}{2 \left(1 + t^2\right)}\nonumber\\
y(t,\theta)&= \frac{1 - p^2 + t^2 + \left(c + b t\right)^2+\sqrt{4\left (1 + t^2\right)p^2 +\left(1 - p^2 + t^2 + \left(c + b t\right)^2\right)^2}\cos\theta}{2 \left(1 + t^2\right)}t\nonumber\\
z(t,\theta)&= \sqrt{p^2 +\frac{\left(1 - p^2 + t^2 + \left(c + b t\right)^2\right)^2}{4 \left(1 + t^2\right)}}\sin\theta, \quad\quad (t, \theta) \in \mathbb R \times [0,2\pi).
\end{align}
After expressing $\cos\theta$ from the first equation
we obtain the following equalities:
$$ 
y= t x, \quad
z^2= -p^2 \left(-1 + x\right) + \left(1 + c^2 + 2 b c t + t^2 \left(1 + b^2 - x\right) - x\right) x.
$$
If we eliminate the parameter $t$ from the above equations, we obtain the following implicit equation of $\mathcal{CS}(\mathbf l,p)$:
\begin{equation}
\label{implicit_parabolic}
x \left(x^2+y^2+z^2\right) -x^2 (c^2 + 1 - p^2) - 2 b c x y - y^2 (b^2 + 1) - p^2 x=0.
\end{equation}
It is clear from eq. (\ref{implicit_parabolic}) that $\mathcal{CS}(\mathbf l,p)$  is an algebraic surface of the order 3 containing the absolute conic and the axis $z$, i.e. $\mathcal{CS}(\mathbf l,p)$ is a parabolic cyclide through $z$.
\hfill {\scriptsize$\square$}

\begin{tm}\label{tm1}
Let $\alpha$ be an $m^{th}$ order algebraic curve that cuts the axis $z$ at $z'$ points, the absolute conic at $a'$ pairs of the absolute points and with the points $P_1$ and $P_2$ as $p'_1$-fold and $p'_2$-fold points, respectively. Then, the following statements hold:
\begin{enumerate} 
\item $\mathcal{CS}(\alpha,p)$ is an algebraic surfaces of the order $3m-(z'+2a'+2p'_1+2p'_2)$.
\item  The  absolute conic is an  $m-(z'+p'_1+p'_2)$-fold curve of $\mathcal{CS}(\alpha,p)$.
\item The axis $z$ is an $(m-2a'+z')$--fold line of $\mathcal{CS}(\alpha,p)$.
\item  The points  $P_1$, $P_2$ are  $2m-(2a'+p_1'+p_2')$--fold points of \,$\mathcal{CS}(\alpha,p)$.
\end{enumerate}
\end{tm}

\noindent {\sc Proof.}

\begin{enumerate}
\item The order of $\mathcal{CS}(\alpha,p)$ equals the number of  intersection points of $\mathcal{CS}(\alpha,p)$ and any straight line $\mathbf l$. Through every such intersection point  passes one circle $c\in \mathcal C(p)$ that cuts $l$ and $\alpha$. Since $c$ must lie on $\mathcal{CS}(\mathbf l,p)$, the number of such circles is equal to the number of  intersection points of $\alpha$ and $\mathcal{CS}(\mathbf l,p)$ that is, according to Lemma \ref{L1}, $3m$. 
 \begin{itemize}
\item[--] If $\alpha$ intersects the axis $z$ at  $z'$ points different from $P_i$ or at a certain number of such points where the sum of their multiplicities is $z'$, the system $\mathcal K(\alpha,p)$  contains the $z'$-ple plane $\omega$. Thus, the number of the intersection points of $\mathbf l$ and $\mathcal{CS}(\alpha,p)$ is reduced by $z'$.
\item[--] If $\alpha$ passes through $a'$ pairs of  absolute points, the  system $\mathcal K(\alpha,p)$ contains $a'$ pairs of  isotropic planes through the axis $z$. Thus, the number of  intersection points of $\mathbf l$ and $\mathcal{CS}(\alpha,p)$ is reduced by $2a'$.
\item[--] If $P_i$ ($i=1,2$) is the $p'_i$-fold point of $\alpha$, the system $\mathcal K(\alpha,p)$  contains $p'_i$ planes through the axis $z$ and the $p'_i$-ple plane $\omega$. Thus, the number of intersection points of $\mathbf l$ and $\mathcal{CS}(\alpha,p)$ is reduced by $2p_i$.
\end{itemize}
These observations lead to the conclusion that the order of $\mathcal {SC}(\alpha,p)$ is $3m-(z'+2a'+2p'_1+2p'_2)$.

\item The plane at infinity cuts the system $\mathcal K(\alpha,p)$ into $m$ straight lines that pass through the points of $\alpha$ at infinity and the absolute conic counted $m$ times.
 \begin{itemize}
\item[--] If $z'\neq 0$, $\mathcal K(\alpha,p)$ contains the $z'$-ple plane $\omega$. Thus, the multiplicity of the absolute conic as the section of $\mathcal {CS}(\alpha,p)$ and the plane $\omega$ is reduced by $z'$.
\item[--]  If $a'\neq 0$,  the number of straight lines at infinity on $\mathcal {CS}(\alpha,p)$ is reduced by $2a'$.
\item[--] If $p'_i\neq 0$, the number of straight lines at infinity on $\mathcal {CS}(\alpha,p)$ as well as the multiplicity of the absolute conic as the section of $\mathcal {CS}(\alpha,p)$ and the plane $\omega$ is reduced by $p'_i$.
\end{itemize}
Therefore,  $\mathcal{SC}(\alpha,p)$ contains the absolute conic counted $m-(z'+p'_1+p'_2)$ times.

\item In the general case (if $\alpha$ does not contain any singular point of $\mathcal C(p)$), every plane $\zeta\in [z]$ cuts $\mathcal {CS}(\alpha,p)$ into $m$ circles and the axis $z$. Thus, $z$ is an $m$-fold line of $\mathcal {CS}(\alpha,p)$.
 \begin{itemize}
\item[--] If $z'\neq 0$, the $z'$-fold plane $\omega$ is excluded from $\mathcal K(\alpha,p)$, but the $z'$-fold axis $z$ (as the part of degenerated circles in $z'$ planes through $z$ and the tangent lines of $\alpha$ at its intersection points with $z$) lie on $\mathcal {CS}(\alpha,p)$.
\item[--]  If $a'\neq 0$, then $2a'$ isotropic planes with the $2a'$-fold axis $z$ are excluded from $\mathcal K(\alpha,p)$.
\item[--] If $p_i'\neq 0$, $p'_i$ planes  with the $p'_i$-fold axis $z$ are excluded from $\mathcal K(\alpha,p)$ but the $p_i'$-fold axis $z$ (as the part of degenerated circles in $p'_i$ planes through $z$ and the tangent lines of $\alpha$ at $P_i$) lie on $\mathcal {CS}(\alpha,p)$.
\end{itemize}
Thus, the axis $z$ is an $(m-2a'+z')$-fold line of $\mathcal {CS}(\alpha,p)$.

\item The order of an algebraic cone is equal to the number of its generators that lie in any plane through its vertex. If a point $A$ is the singular point of a surface $\mathcal S$ and lies  on its $n$-fold straight line $a$, then the tangent cone of $\mathcal S$ at $A$ contains $a$ as an $n$-fold generator.

\noindent Since every plane $\zeta\in [z]$ cuts $\mathcal {CS}(\alpha,p)$ into the $(m-2a'+z')$-fold line $z$, and $m-(z'+p'_1+p'_2)$ circles through $P_i$, this plane cuts the tangent cone of $\mathcal {CS}(\alpha,p)$ at $P_i$ into the $(m-2a'+z')$-fold line $z$ and $m-(z'+p'_1+p'_2)$ tangents of the intersection circles. Thus, the order of the tangent cone of $\mathcal {CS}(\alpha,p)$ at $P_i$ is $2m-(2a'+p'_1+p'_2)$.\hfill$\square$
\end{enumerate}

\begin{cor}
If $\alpha$ is a planar algebraic curve that lies in any plane $\zeta (\varphi)$, then $\mathcal {CS}(\alpha,p)$ is identical to the plane $\zeta (\varphi)$ counted $2m-(2a'+p'_1+p'_2)=m-2a'+z'$ times. 
\end{cor}

\noindent {\sc Proof.}
The conclusion follows from the facts that for every point $A \in \alpha$, the circle $c^A(p)\in\mathcal C(p)$ lies in the plane $\zeta (\varphi)$, and  $m=z'+p'_1+p'_2$. 
\hfill$\square$

\begin{cor}
If $\alpha$ is an algebraic curve that lies  on any torus $\tau (r)$, then $\mathcal {CS}(\alpha,p)$ is identical to the torus $\tau (r)$ counted $\displaystyle \frac{m-p'_1-p'_2}{2}$ times.
\end{cor}

\noindent {\sc Proof.}
For every point $A \in \alpha$, the circle $c^A(p)\in\mathcal C(p)$ also lies on the torus $\tau (r)$. Since the torus is a surface of the fourth order having  double points in $P_1, P_2$, the numbers $m, p_1', p'_2$ are even and $z'=0, 2a'=m$.\\
Particularly, for $q>0$, $\tau(p)$ is a double sphere $S(O,p)$ through the points $P_1$ and $P_2$. Therefore,  if $\alpha$ is an algebraic curve that lies on the sphere $S(O,p)$, then $\mathcal {CS}(\alpha,p)$ is identical to the sphere $S(O,p)$ counted $m-p'_1-p'_2$ times.
\hfill$\square$

\bigskip

\noindent {\bf Singular lines} and {\bf singular points} of $\mathcal {CS}(\alpha,p)$\\
Some of the singular lines and singular points of the circular surface $\mathcal {CS}(\alpha,p)$ have already been mentioned in this paper. Besides them, $\mathcal {CS}(\alpha,p)$ can have some further singularities. Here we give the complete list of singular lines and points of $\mathcal {CS}(\alpha,p)$:
\begin{itemize}
\item[-] The line $z$ and the absolute conic are in the general case the singular lines of $\mathcal {CS}(\alpha,p)$. All points lying on $z$, the points $P_1, P_2$ and the absolute points are the singular points of $\mathcal {CS}(\alpha,p)$.

\item[-]
 If different points of the curve $\alpha$ determine the same circle $c\in\mathcal C(p)$, this circle is a singular circle of  $\mathcal {CS}(\alpha,p)$ and its multiplicity equals the number of  intersection points of $\alpha$ and $c$. Through the singular point of $\alpha$, the singular circle of $\mathcal {CS}(\alpha,p)$ passes.

\item[-] Let $A\in \mathcal C(p)$ be a point that does not lie on any singular line of  $\mathcal C(p)$. The necessary condition for the singularity of $A$ is that the circle $c_A$ has the singular point at $A$. It is possible only  if $\mathcal C(p)$ is a hyperbolic congruence ($q<0$) and $A$ lies on the circle $c (0)$ that is given by the equations $x^2+y^2=-p^2$, $z=0$. Therefore, $A$ is the singular point of $\mathcal C(p)$  only if it is the intersection of $\alpha$ and $c(0)$. 
\end{itemize}

\bigskip

\noindent {\bf Remark}\\
In the case of an elliptic congruence $\mathcal C(p)$, $q>0$, a circular surface $\mathcal {CS}(p,\alpha)$ can be obtained by a spherical inversion, i.e. by  a quadratic inversion with respect to a sphere and its center as a pole. Namely, if $i$ is a spherical inversion defined by the sphere $S(P_1,\overline{P_1P_2})$, then $\mathcal {CS}(p,\alpha)$ is the inverse image of $\mathcal B(P_2,i(\alpha))$ that is the cone with the vertex $P_2$ and the base curve $i(\alpha)$ (see Fig. 8).

\begin{center}
\includegraphics[scale=0.8]{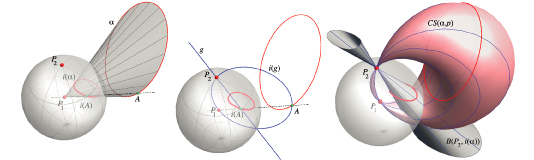}
\vspace{-0.5cm}
\begin{multicols}{3}
a

b

c
\end{multicols}
\end{center}
\vspace{-0.2cm}
\noindent Figure 8: \parbox[t]{120mm}{The figure illustrates the following: if $i$ is a spherical inversion defined by the sphere $\mathcal S(P_1,\overline{P_1P_2})$ and if $\alpha$ is a curve, every generatrix $g$ of the cone $\mathcal B(P_2,i(\alpha))$ is transformed to a circle  passing through the points $P_1$ and $P_2$ that intersect a curve $\alpha$, i.e. $i(\mathcal B(P_2,i(\alpha))=\mathcal {CS}(p,\alpha)$. The statement is also valid for the changed positions of  $P_1$ and $P_2$.}

\medskip

\noindent All properties of $\mathcal {CS}(p,\alpha)$  listed in Theorem \ref{tm1} follow as the properties of $i(\mathcal B(P_1,i(\alpha))$ from the  properties of the spherical inversion $i$  (see \cite{Nice}).

\subsection{Algebraic equation of $\mathcal {CS}(\alpha,p)$}

Let the poynomial $x^2+y^2+z^2$ be denoted by $A$. Let $F^n=F^n(x,y,z)$ and $H^n=H^n(x,y,z)$ denote the $n^{\rm {th}}$ degree polynomials in $x$, $y$ and $z$ with real coefficients, and let $H^n$ be a homogeneous polynomial. Then, by analogy with the equation of the $s$-circular planar curve given in \cite{wiki_circular}, the implicit equation of a surface that contains the absolute conic as an $s$-fold curve is the following:

\begin{equation}
\label{implicit_general}
F^n=\sum_{i=0}^{s-1} H^{n-2s+i}A^{s-i}+\sum_{j=0}^{n-s}H^j=0, \quad H^{n-2s}\neq 0,\,\,n\geqslant 2s.
\end{equation}
\noindent Further, according to the theorem given in \cite[p.\,251]{Harris}, if an $n$-th order surface in $\mathbb E^3$ which passes through the origin is given by the equation
\begin{equation}
\label{tangent_cone}
F^n=H^k+H^{k+1}+\cdots + H^n= 0,  \quad 1\leqslant k \leqslant n,
\end{equation}
 then the tangent cone at the origin is given by $H^k(x, y, z) = 0$.

\medskip

\noindent According to Theorem \ref{tm1}, for the circular surface $\mathcal {CS}(\alpha,p)$, the numbers $n$ and $s$ from eq. (\ref{implicit_general}) are equal to $3m-(z'+2a'+2p'_1+2p'_2)$ and $m-(z'+p'_1+p'_2)$, respectively, and the multiplicity of the origin is at least $m+z'-2a'$.  Therefore, the implicit equation of $\mathcal {CS}(\alpha,p)$ takes the following form:
\begin{equation}
\label{implicit_circular}
\sum_{i=0}^{m-(z'+p'_1+p'_2+1)} H^{m+z'-2a+i}A^{m-(z'+p'_1+p'_2+i)}+\sum_{j=m-2a'+z'}^{2m-(2a'+p'_1+p'_2)}H^j=0.
\end{equation}

\section{Examples}

In this section we give some examples of the circular surfaces  $\mathcal {CS}(\alpha,p)$ when $\alpha$ is a line, a conic, a twisted cubic or a cyclic harmonic curve. For a particular circular surface  we obtain the implicit equation (eliminating the parameters $t$ and $\theta$ from the parametric equations (\ref{par.surface1})), check the properties given in Theorem \ref{tm1} and visualize the shape of the surface. For  computing and plotting, we use  the program {\it Mathematica}.

\subsection{ $\alpha$ is a straight line}
It was shown in Lemma \ref{L1} that if $\alpha$ is a straight line given by $\mathbf l(t)=(1,t,bt+c)$, $t,b,c\in\mathbb R$, the implicit equation of $\mathcal {CS}(\mathbf l,p)$ is given by eq. (\ref{implicit_parabolic}), i.e. $\mathcal {CS}(\mathbf l,p)$ is a parabolic cyclide.
 According to the theorem related to (\ref{tangent_cone}),  the tangent cone of  $\mathcal {CS}(\mathbf l,p)$ at $P_1$ ($P_2$), in the coordinate system with the origin at $P_1$ ($P_2$) is given by the  equation $(1+c^2-p^2)x^2+2 b c x y+y^2+b^2 y^2 \mp 2 p x z=0$.

\begin{center}
\includegraphics[scale=0.8]{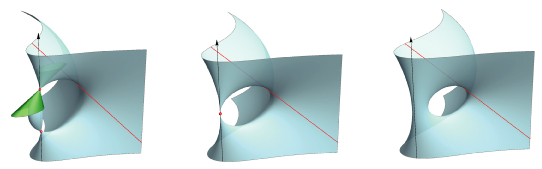}
\vspace{-0.5cm}
\begin{multicols}{3}
a

b

c
\end{multicols}
\vspace{-0.2cm}
\noindent Figure 9: \parbox[t]{120mm}{$\mathcal {CS}(\mathbf l,p)$ for $\mathbf l(t)=(1, t,t+2)$ and $p=1,0,i$ are shown in  a, b and c, respectively. In a the tangent cone at the point $P_2$ is also shown.
}
\end{center}

\subsection{ $\alpha$ is a conic}

If $\alpha$ is a conic in the general position to the singular points of $\mathcal C(p)$, the circular surface $\mathcal {CS}(\alpha,p)$ is a sextic with the double axis $z$ and the double absolute conic. The points $P_i$ are the quadruple points of $\mathcal {CS}(\alpha,p)$.

\noindent The directing curve for the circular surfaces in Figure 10 is a rectangular hyperbola  $\mathbf h_1(t)=(t,\frac{1}{t},0)$,\, $ t\in \mathbb R.$ The implicit equation of $\mathcal {CS}(\mathbf h_1,p)$  is the following:
\begin{equation}
\label{h1}
x y A^2-(x^4 + 2 x^2 y^2 + p^4 x^2 y^2 + y^4 + 2 p^2 x y z^2)+p^4 x y=0.
\end{equation}
\noindent Every point $Z\in z$, $Z\neq P_i$, is the biplanar point of $\mathcal {CS}(\mathbf h_1,p)$ with real tangent planes $x=0$, $y=0$. The tangent cone of $\mathcal {CS}(\mathbf h_1,p)$ at $P_i$ in the coordinate system with an origin at $P_i$ is given by   $
x^4 - 2 p^2 x^3 y + 2 x^2 y^2 + p^4 x^2 y^2 - 2 p^2 x y^3 + y^4 - 
 4 p^2 x y z^2=0.
$ For $q>0$ this cone is a proper real quartic cone with the double generatrix  $z$ (see Fig. 10a).

\begin{center}
\includegraphics[scale=0.92]{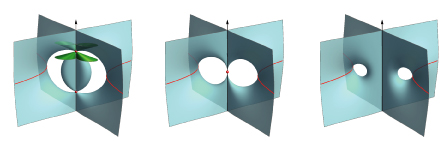}
\vspace{-0.5cm}
\begin{multicols}{3}
a

b

c
\end{multicols}
\end{center}
Figure 10:  \parbox[t]{120mm}{$\mathcal {CS}(\mathbf h_1,p)$ for $p=1,0,i$ are shown in a, b, c, respectively.
}
\bigskip

\noindent If  $\alpha$ is a conic that cuts the axis $z$ at one point $Z\neq P_i$, the circular surface $\mathcal {CS}(p,\alpha)$ is a quintic  with a triple axis $z$,  quadruple points $P_i$ and a simple absolute conic. 

\noindent The directing curve for the circular surfaces in Figure 11 is a rectangular hyperbola $\mathbf h_2$ that lies in the plane $x=1$ and cuts the axis $z$ at the point at infinity, $\mathbf h_2(t)=(1,t,\frac{1}{t})$,\, $ t\in\mathbb R.$
 The implicit equation of $\mathcal {CS}(\mathbf h_2,p)$ is the following:
\begin{equation}
x y^2 A-(x^4 + x^2 y^2 - p^2 x^2 y^2 + y^4)-p^2 x y^2=0.
\end{equation}
The plane at infinity intersects $\mathcal {CS}(\mathbf h_2,p)$ at the absolute conic, the line in the plane $y=0$ counted twice and the simple line in the plane $x=0$. The axis  $z$ is the triple line of $\mathcal {CS}(\mathbf h_2,p)$ with the splitting tangent cone at every point $(0,0,z_0)$, $z_0\neq \pm p$. This tangent cone splits into the plane $y=0$ counted twice and the plane $x=0$. In the coordinate system with an origin in $P_i$, the tangent cone of $\mathcal {CS}(\mathbf h_2,p)$ at $P_i$ is given by the equation $x^4 + x^2 y^2 - p^2 x^2 y^2 + y^4 - 2 p x y^2 z=0.$
This equation presents a proper real quartic cone (see Fig. 11a) with a triple generatrix  $z$, for $q>0$, and  two pairs  of imaginary planes through $z$, for $p=0$.

\begin{center}
\includegraphics[scale=0.8]{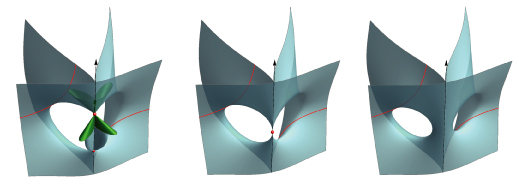}
\vspace{-0.5cm}
\begin{multicols}{3}
a

b

c
\end{multicols}
\end{center}
\vspace{-0.2cm}
Figure 11:  \parbox[t]{120mm}{$\mathcal {CS}(\mathbf h_2,p)$ for $p=1,0,i$ are shown in a, b, c, respectively.
}

\bigskip

\noindent In general, if $\alpha$ is a circle, $\mathcal {CS}(\alpha,p)$ is a cyclide, i.e. a quartic passing two times through the absolute conic.  $P_i$ are the double points of $\mathcal {CS}(\alpha,p)$. Particularly, if a circle lies on any torus $\tau (r)$, its circular surface is $\tau (r)$. If a circle cuts the axis $z$ at any point $Z\neq P_i$, the circular surface is a parabolic cyclide.  If a circle passes through any point $P_i$, its circular surface is a sphere through $P_1$ and $P_2$. If a circle lies on any sphere through $P_1$ and $P_2$, its circular surface is this sphere, or a part of it, counted twice. Therefore, for any circle that does not lie in a plane $\zeta (\varphi )$, there exists a number $p$ such that the circular surface of this circle is a double sphere or a double part of a sphere. Namely, for any circle there exists a sphere that passes through it and with the center in the plane $z=0$. The axis $z$ cuts this sphere at the points $(0,0,\pm p)$.

\bigskip

\noindent For $q<0$ and a conic $\alpha$,  $\mathcal {CS}(\alpha,p)$ can possess 1, 2, 3 or 4 double points that do not lie on its singular lines. Some examples are shown in Figure 12.

\begin{center}
\includegraphics[scale=0.74]{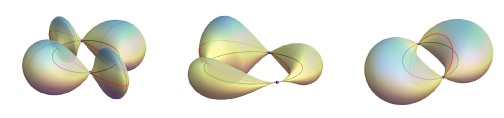}
\vspace{-0.5cm}
\begin{multicols}{3}
a

b

c
\end{multicols}
\end{center}
\vspace{-0.2cm}
Figure 12:  \parbox[t]{120mm}{$\mathcal {CS}(\alpha,p)$ with 4, 3 and 2 double points out of the axis $z$ are shown in  a, b and c, respectively. In  a,  $\alpha(t)=(\frac{1}{2}\cos t,2\sin t,0)$ and $p=i$. In  b, $\alpha(t)=(\frac{3}{2}\cos t,3\sin t-1,0)$ and $p=2i$. In  c, $\alpha=\mathbf c_1(t)=(\frac{\sqrt{2}}{2}\cos t, \sin t, \frac{\sqrt{2}}{2}\cos t),$ and $p=i$.}

\subsection{$\alpha$ is a twisted cubic}

In the general case, if $\alpha$ is a twisted cubic, $\mathcal {CS}(\alpha,p)$ is a surface of the order 9  with 3-ple lines at the absolute conic and axis $z$, and  with  6-fold points at $P_i$. If $\alpha$ passes through the singular points of $\mathcal C(p)$, there are many cases when the order of $\mathcal {CS}(\alpha,p)$ decreases.  The minimal order of such surfaces is 3 in the case when $\alpha$ is a twisted cubic circle that passes through the points $P_i$.

\noindent The directing curve for the circular surfaces in Figure 13 is a twisted cubic circle (\cite[pp. 69-76]{Fladt})  given by the  parametrization $\mathbf t(t)=(\frac{t}{1 + t^2},\frac{ t + t^3}{1 + t^2}, \frac{t^2}{1 + t^2})$,\, $ t\in\mathbb R$.
By eliminating the parameters $t$ and $\varphi$ from eqs. (\ref{par.surface1}) for $\alpha=\mathbf t$, we obtain the following algebraic equation:
\begin{align}
\label{twisted_circle}
&x (x - y) A^2\nonumber\\
+&(x^2 - y^2)^2 -
 p^2 (2 x^4 - 4 x^3 y + 2 x^2 y^2 + 2 x^2 z^2 - 2 x y z^2)+ p^4 x^2 y^2 \nonumber\\
+&p^4 x (x - y)=0.
\end{align}

\noindent For $p\neq 0$ this equation presents $\mathcal {CS}(\mathbf t,p)$ that is a 6$^\text{th}$ order surface, with a double absolute conic and the axis $z$, and with quadruple points at $P_i$.
The quadruple tangent cones of $\mathcal {CS}(\mathbf t,p)$ at real points $P_i$, in the coordinate systems with the origin at $P_i$, are given by $
(x^2 - y^2)^2 + p^4 (x^2 y^2) + 2 p^2 x (x - y) (x y + y^2 + 2 z^2) =0.
$

\noindent For $p=0$, i.e. when $\mathbf t$ cuts the axis $z$ at the point $P$, eq. (\ref{twisted_circle}) presents the splitting 6$^\text{th}$ order surface that splits into the plane $x=y$ and the circular surface $\mathcal {CS}(\mathbf t,0)$ that is a 5$^\text{th}$ order surface with a double absolute conic, a torsal line through the axis $z$ (the tangent plane is $x=0$), and a triple point at the origin. The algebraic equation of $\mathcal {CS}(\mathbf t,0)$ is the following:
\begin{equation}
x A^2+(x - y) (x + y)^2=0.
\end{equation}

\begin{center}
\includegraphics[scale=0.73]{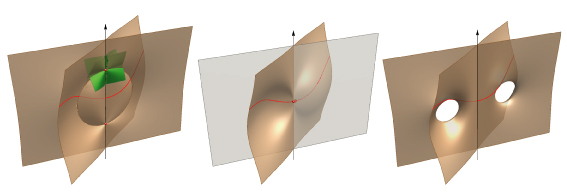}
\vspace{-0.7cm}
\begin{multicols}{3}
a

b

c
\end{multicols}
\end{center}
\vspace{-0.2cm}
Figure 13:  \parbox[t]{120mm}{$\mathcal {CS}(\mathbf t,p)$ for $p=1,0,i$ are shown in  a, b, c, respectively.
}

\subsection{$\alpha$ is a cyclic-harmonic}

\noindent The {\it rose surfaces} studied in \cite{roses} are circular surfaces $\mathcal {CS}(\alpha,p)$ where $\alpha$  is a rose (rhodonea curve) given by the polar equation $\rho=\cos \frac{n}{d}\varphi$, and  $\mathcal C(p)$ is an elliptic or parabolic congruence. The rose lies in the plane perpendicular to $z$ having the directing point $P_i$ as its multiple point.

\noindent If we extend $\alpha$ to all cyclic-harmonic curves $\rho=\cos \frac{n}{d}\varphi+k$, $k\in \mathbb R^+\cup \{0\}$ (see \cite{Hilton/clanak},\cite {Moritz}),  and include a  hyperbolic congruence  $\mathcal C(p)$, numerous forms of a new class of surfaces with nice visualizations and interesting algebraic properties can be obtained. The authors suppose that it is worth studying a complete class of such surfaces in a further work. Here, we show only one example for one of the simplest  cyclic-harmonic curves.

\noindent The directing curve $\mathbf {ch}$ for the circular surfaces in Figure 14 is a curtate cyclic-harmonic curve in the plane $z\!=\!0$ given by the polar equation $\rho=\cos 3t+k$, $k>1$, which possesses an isolated 6-fold singularity at the origin.  According to the properties of cyclic harmonic curves and Theorem \ref{tm1}, they are 14$^\text {th}$ order surfaces with a double absolute conic, 10-fold axis $z$ and 12-fold imaginary points $P_i$.

\begin{center}
\includegraphics[scale=0.95]{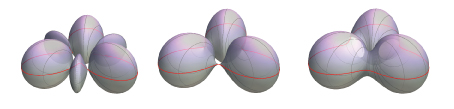}
\vspace{-0.7cm}
\begin{multicols}{3}
a

b

c
\end{multicols}
\end{center}

Figure 14:  \parbox[t]{120mm}{ The surfaces $\mathcal {CS}(\mathbf {ch},i)$ for $k=\frac{3}{2}, 2, 3$ are shown in  a, b and c, respectively. 
}


\begin{thebibliography}{99}

\bibitem{Eisenhart} {\sc L. P. Eisenhart}, Congruences of Curves.  {\it Trans. Amer. Math. Soc.} {\bf 4}, 470--488 (1907).

\bibitem{Fladt} {\sc K. Fladt, A. Baur}, {Analytische Geometrie spezieller Fl\"achen und Raumkurven}. Friedr. Vieweg \& Sohn, Braunschveig, (1975) 

\bibitem{roses} {\sc S. Gorjanc}, Rose Surfaces and their Visualizations. {\it Journal for Geometry and Graphics}, {\bf 13} (1),  59--67 (2010).

\bibitem{Harris} {\sc J. Harris}, {\it Algebraic Geometry}. Springer,
  New York, 1995. 

\bibitem{Hilton/clanak} {\sc H. Hilton}, On Cyclic-Harmonic Curves. {\it The Annals of Mathematics}, Second Series, {\bf 24} (3), 209--212, (1923). \\
{\small  \url{http://www.jstor.org/stable/1967850}}

\bibitem{Izumiya} {\sc 
S. Izumiya, K. Saji, N. Takeuchi}, Circular Surfaces. {\it
Advances in Geometry} {\bf 7} (2),  295--313 (2007).

\bibitem{Moritz} {\sc R. E. Moritz}  On the Construction of Certain Curves Given in Polar Coordinates. {\it American Mathematical Monthly} {\bf 24} (5), 213-220 (1917).

\bibitem{Nice}
{\sc V. Ni\v{c}e},  Krivulje i plohe 3. i 4. reda nastale pomo\'{c}u kvadratne inverzije. {\it Rad HAZU} {\bf 278} (86), 153--194 (1945).


\bibitem{Salmon2} {\sc G. Salmon}, {\it A Treatise on the Analytic
    Geometry of Three Dimensions},  Vol.II., Chelsea Publishing Company,
    New York, (reprint), 1965.

\bibitem{wiki_Apollonian}
{\small \url{http://en.wikipedia.org/wiki/Apollonian_circles}}

\bibitem{wiki_circular}
{\small \url{http://en.wikipedia.org/wiki/Circular_algebraic_curve}}


\end{thebibliography}
\end{document}